\documentclass[a4paper]{amsart}

\usepackage{amssymb}
\usepackage{amsmath}
\usepackage{amsthm}
\usepackage{amsfonts}
\usepackage{geometry}
\usepackage{graphicx}
\usepackage{subfig}
\usepackage{cases}
\usepackage{txfonts}
\usepackage{todonotes}
\usepackage{braket}
\ifpdf
  \usepackage[pdftex]{hyperref}
  \hypersetup{plainpages=false,unicode=true,pdffitwindow=true}
\else
  \usepackage[dvips]{hyperref}
\fi

\theoremstyle{plain}
\newtheorem{theorem}{Theorem}
\newtheorem{corollary}[theorem]{Corollary}

\theoremstyle{remark}
\newtheorem{remark}[theorem]{Remark}

\newcommand{\R}{\mathbb{R}}

\newcommand{\f}[2]{\frac{#1}{#2}}
\newcommand{\pa}{\partial}
\newcommand{\half}{\f{1}{2}}

\newcommand{\dt}{\Delta t}

\newcommand{\norm}[1]{\left\lVert#1\right\rVert} 

\begin{document}

\title{On the probability density function of baskets}

\author{Christian Bayer}
\address{Weierstrass Institute}
\email{christian.bayer@wias-berlin.de}

\author{Peter K.~Friz}
\address{Weierstrass Institute \and TU Berlin}
\email{friz@math.tu-berlin.de}

\author{Peter Laurence}
\address{Universit\`a di Roma \and Courant Institute}
\email{laurencepm@yahoo.com}

\begin{abstract}
  The state price density of a basket, even under uncorrelated Black--Scholes
  dynamics, does not allow for a closed form density.  (This may be rephrased
  as statement on the sum of lognormals and is especially annoying for such
  are used most frequently in Financial and Actuarial Mathematics.) In this
  note we discuss short time and small volatility expansions,
  respectively. The method works for general multi-factor models with
  correlations and leads to the analysis of a system of ordinary (Hamiltonian)
  differential equations. Surprisingly perhaps, even in two asset
  Black--Scholes situation (with its flat geometry), the expansion can
  degenerate at a critical (basket) strike level; a phenomena which seems to
  have gone unnoticed in the literature to date.  Explicit computations relate
  this to a phase transition from a unique to more than one ``most-likely''
  paths (along which the diffusion, if suitably conditioned, concentrates in
  the afore-mentioned regimes). This also provides a (quantifiable)
  understanding of how precisely a presently out-of-money basket option may
  still end up in-the-money.
\end{abstract}


\maketitle

\section{Introduction}
\label{sec:introduction}

As is well known, the sum of independent log-normal variable does not admit
a closed-form density. And yet, there are countless applications in Finance
and Actuarial Mathematics where such sums play a crucial role, consider for
instance the law of a Black--Scholes basket $B$ at time $T$, i.e. the
weighted average of $d$ geometric Brownian motions.

As a consequence, there is a natural interest in approximations and
expansions, see e.g.~\cite{Du} and the references therein. This article
contains a detailed investigation in small volatility and short time
regimes. Forthcoming work of A.~Gulisashvili and P.~Tankov \cite{GuTa13} deals
with tail asymptotics. Our methods are not restricted to the geometric
Brownian motion case: in principle, each Black--Scholes component could be
replaced by the asset price in a stochastic volatility model, such as the the
Stein--Stein model \cite{Stein91}, with full correlation between all assets
and their volatilities. In the end, explicit solutions only depend on the
analytical tractability of a system of ordinary differential equations. If
such tractability is not given, one can still proceed with numerical ODE
solvers.

As a matter of fact, our aim here is not to push the generality in which our
methods work: one can and should expect involved answers in complicated
models. Rather, our main -- and somewhat surprising -- insight is that
\textit{unexpected phenomena} are already present in the \textit{simplest
  possible setting}: to this end, our first focus will be on the case of $d=2$
independent Black--Scholes assets, without drift and correlation, with unit
spot and unit volatility). To be more specific, if $C_{B}$ denotes the fair
value of an (out-of-money) call option on the basket $B$ struck at $K$, one
naturally expects, for a small maturity $T$,
\begin{equation*}
  \frac{\partial ^{2}}{\partial K^{2}}C_{B}\left(
    K,T\right) \sim \text{\textrm{(const)}} \exp \left( -\frac{\Lambda \left(
        K\right)}{T}\right) \frac{1}{\sqrt{T}}. 
\end{equation*}
And yet, while true for \textit{most} strikes, it fails for $K=K^{\ast }$; in
fact,
\begin{equation*}
\left\{ \frac{\partial ^{2}}{\partial K^{2}}C_{B}\left( K,T\right) \right\}
_{K=K^{\ast }}\sim \text{\textrm{(const)}}\exp \left( -\frac{\Lambda \left(
K^{\ast }\right) }{T}\right) \frac{1}{T^{3/4}}.
\end{equation*}
To the best of our knowledge, and despite the seeming triviality of the
situation (two independent Black--Scholes assets!), the existence of a
``special'' strike level $K^{\ast }$, at which the value of a basket option
(here: butterfly spread\footnote{Extensions to spreads and vanilla options are
  possible and will be discussed elsewhere.}) has a ``special'' decay
behavior, as maturity approaches $0$, seems to be new. There are different
proofs of this fact; the most elementary argument -- based on the analysis of
a convolution integral -- is given in
Section~\ref{sec:comp-based-saddle}. However, this approach -- while telling
us {\it what} happens -- does not tell us {\it how} it happens.

The main contribution of this note is precisely a good understanding of the
latter.  In fact, there is clear picture that comes with $K^{\ast }$.  For
$K<K^{\ast }$ and conditional on the option to expire on the money, there is a
unique ``most likely'' path around which the underlying asset price process
will concentrate as maturity approaches $0$. For $K>K^{\ast }$, however, this
ceases to be true: there will be two distinct (here: equally likely) paths
around which concentration occurs.  What underlies this interpretation is that
large deviation theory
not only characterizes the probability of unlikely events (such as expiration
in-the-money, if presently out-of-the-money, as time to maturity goes to zero)
but also the mechanism via which these events can occur.  Such understanding
was already crucial in previous works on baskets aiming at quantification of
basket (implied vol) skew relative to its components, starting with
\cite{A1,A2}. As a matter of fact, the analysis in these papers relied on the
statement that ``generically there is a unique arrival point [of a unique
energy minimizing path] on the (basket-strike) arrival manifold''. The
situation, however, even in the Black--Scholes model, is more involved. And
indeed, we shall establish existence of a critical strike $K^{\ast }$, at
which one sees the phase-transition from one to two energy minimizing, ``most
likely'', paths.\footnote{It can be shown that, sufficiently close to the
  arrival manifold, there is in fact a unique energy minimizing paths. The
  (near-the-money) analysis of \cite{A1,A2} is then justified.} And this
information will have meaning to traders (as long as they believe in a
diffusion model as maturity approaches $0$, which may or may not be a good
idea \ldots) as it tells them the possible scenarios in which an out-of-the
money basket option may still expire in the money.

Let us conclude this introduction with a few technical notes. We view the
evolution of the basket price -- even in the Black-Scholes model -- as a
stochastic volatility evolution model; by which we mean $dB_t/B_t =
\sigma(t,\omega) dW_t$ (as opposed to a local vol evolution where
$\sigma=\sigma(t,B_t)$). This should explain why the methods developed in
Part I of \cite{DFJVpartI,DFJVpartII} for the analysis of stochastic
volatility models (then used in Part II, \cite{DFJVpartII}, to solve the
concrete smile problem (shape of the wings) for the correlated Stein--Stein
model), are also adequate for the analysis of baskets. In a sense, the present
note may well be viewed as Part III in this sequence of papers.

{\it Acknowledgment:} Martin Forde kindly informed us about some misleading
formulations in a previous version.
 P.K.F. has received partial funding from the European
Research Council under the European Union's Seventh Framework Program
(FP7/2007-2013) / ERC grant agreement nr. 258237.

\section{Computations based on saddle-point method}
\label{sec:comp-based-saddle}

In terms of a standard  $d$-dimensional Wiener process $\left(
W^{1},\dots ,W^{d}\right) $,%
\begin{equation*}
B_{T}=\sum_{i=1}^{d}S_{0}^{i}\exp \left( \mu ^{i}T+\sigma
^{i}W_{T}^{i}\right) .
\end{equation*}%
Write $f=f_{T}\left( K\right) $ for the probability density function of $
B_{T}$; i.e. for $\mathbb{P}\left[ B_{T}\in [ K,K+dK]\right] /dK$. Of
course, it is given by some $\left( d-1\right) $-dimensional convolution
integral, explicit asymptotic expansions are -- in principle -- possible with
the saddle point method. It will be enough for our purposes to illustrate
the method in the afore-mentioned simplest possible setting:%
\begin{equation*}
d=2,\,S_{0}^{1}=S_{0}^{2}=1,\,\mu ^{1}=\mu ^{2}=0,\,\sigma ^{1}=\sigma
^{2}=1.
\end{equation*}
In other words, $B_{T}=\exp \left( W_{T}^{1}\right) +\exp \left(
W_{T}^{2}\right)$. 
We claim that for some constant  $c_0 = c_0(K) > 0 $
\begin{subequations}
  \label{eq:density-expansion}
  \begin{numcases}{f\left( K\right) =}
      \exp \left( -\frac{\Lambda \left( K\right) }{T}\right)
      \frac{1}{\sqrt{T}} \left( c_0+O\left( T\right) \right), & when $K
      \neq K^{\ast }$, \label{eq:density-expansion-a}\\ 
      \exp \left( -\frac{\Lambda \left( K^{\ast }\right) }{T}\right)
      \frac{1}{T^{3/4}}\left( c_0+O\left( T\right) \right), &  when $K =
      K^{\ast }$, \label{eq:density-expansion-b}
  \end{numcases}
\end{subequations}
with
\begin{equation*}
K^{\ast }=2e\approx 5.43656
\end{equation*}
and
\begin{equation*}
  \Lambda(K) = \inf \{ h_K(x) \, | \, x \in [0,K]\} 
\end{equation*}
with
\begin{equation}
  \label{eq:h-function}
  h_K(x) := \frac{1}{2} \left( \left(\log x \right)^2 + \left(\log(K-x)\right)^2 \right).
\end{equation}
Note that for $K \le K^\ast$ we can explicitly solve this minimization problem
and obtain $\Lambda(K) = \log(K/2)^2$ with corresponding minimizer $x^\ast =
K/2$, corresponding to the single local extremum of $h_K$. For $K > K^\ast$,
we have two global minima, which cannot be given in closed form, and hence
$\Lambda(K)$ can only be computed numerically.
\begin{figure}[!htb]
  \centering
  \subfloat[$h_K$ for $K<K^\ast$]{
    \includegraphics[width=0.28\textwidth]{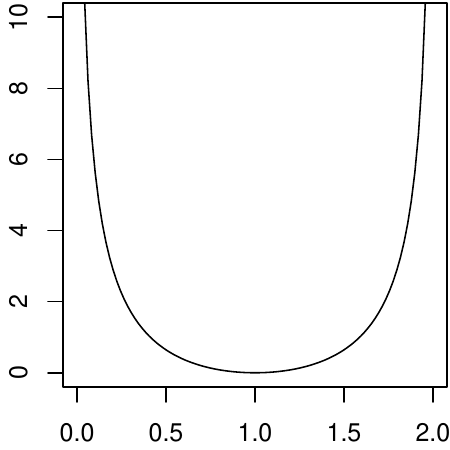}
  }
  \qquad
  \subfloat[$h_K$ for $K=K^\ast$]{
    \includegraphics[width=0.28\textwidth]{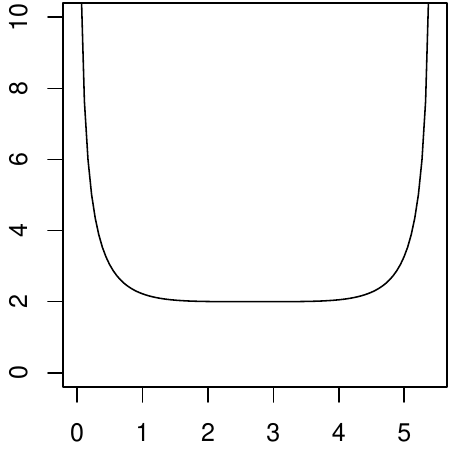}
  }
  \qquad
  \subfloat[$h_K$ for $K>K^\ast$]{
    \includegraphics[width=0.28\textwidth]{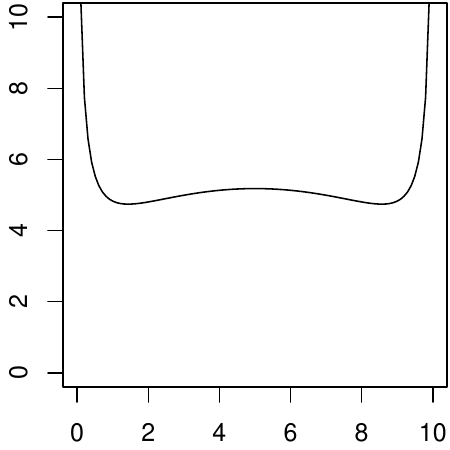}
  }
  \caption{Plot of $h_K$ for different choices of $K$.\\
    (A) For $K < K^\ast$ there is a unique global minimum at $x^\ast = K/2$ which
    is non-degenerate in the sense that $h^{\prime\prime}(x^\ast) > 0.$ \\ 
    (B) For $K = K^\ast$ there is a unique global minimum at $x = K/2$ which
    is degenerate in the sense that $h^{\prime\prime}(x^\ast) = 0.$\\
    (C) For $K > K^\ast$, $x = K/2$ gives a local maximum. There are two
    symmetric global minimizers, which are not given in closed form.}
  \label{fig:hK-plot}
\end{figure}

The stock price $S^i_T$ has a log-normal distribution with parameters $\mu^i =
0$ and $\xi^i = \sigma^i \sqrt{T} = \sqrt{T}$, where the density of the
log-normal distribution is given by
\begin{equation}
  \label{eq:lognormal-dens}
  f_{\mu,\xi}(x) = \frac{1}{\sqrt{2\pi}\xi x} \exp\left( - \frac{\left(\log x
        - \mu  \right)^2}{2 \xi^2} \right).
\end{equation}
Obviously, the density of the sum of these two independent log-normal random
variables satisfies
\begin{equation}
  \label{eq:convolution}
  f(K) = \int_0^K f_{\mu^1,\xi^1}(K-x) f_{\mu^2,\xi^2}(x) dx.
\end{equation}
Using our special parameters, the integrand is of the form
\begin{equation*}
  f_{\mu^1,\xi^1}(K-x) f_{\mu^2,\xi^2}(x) = \frac{1}{2 \pi T
    x(K-x)} \exp\left(-\frac{h_K(x)}{2T} \right).
\end{equation*}
In order to apply the Laplace approximation to~\eqref{eq:convolution}, we
compute the minimizer for $h_K$, which is found by the first order condition
\begin{equation}  \label{eq:LaplaceFirstOrderCond}
  h_K'(x) = 0 \iff \frac{\log x}{x} - \frac{\log(K-x)}{K-x} = 0.
\end{equation}
Clearly, this equation is solved by choosing $x^\ast = K/2$ --- which is the
unique global minimizer iff $K \le 2e$ and a local maximizer otherwise, in which 
case we have two global minima $x_1^\ast < K/2 < x_2^\ast$.
Assuming $K \le 2e$, we can check  degeneracy of that minimum directly by computing
\begin{equation*}
  h_K''(x^\ast) = h_K''(K/2) = 16 \frac{1 - \log(K/2)}{K^2}.
\end{equation*}
and
\begin{equation}
  \label{eq:degeneracy-saddle-point}
  h_K''(x^\ast) = 0 \iff K = 2 e.
\end{equation}
With more work one can see that also the global minima $x_1^\ast, x_2^\ast$, in
the case $K>2e$, are non-degenerate. Hence, whenever $K \ne 2 e$ a standard Laplace method
leads to the expansion \eqref{eq:density-expansion-a}. In the remainder of this section, we consider
the degenerate case and establish \eqref{eq:density-expansion-b}.

Choosing $K = 2 e$ and, correspondingly, $x^\ast = e$, we obtain the Taylor
expansion $h_K(x) = h_K(x^\ast) + \frac{h_K^{(4)}(x^\ast)}{24} (x-x^\ast)^4 +
\mathcal{O}((x-x^\ast)^5)$, with $h_K(x^\ast) = 2$ and $h_K^{(4)}(x^\ast) = 20
e^{-4}$, we obtain the Laplace approximation
\begin{align*}
  f(K) &= \int_0^K \frac{1}{2\pi T (K-x)x} \exp\left(-\frac{h_K(x)}{2T}
  \right) dx\\
  &= \frac{1}{2\pi T e^{2}} \int_0^K \exp\left(-\frac{1}{T}\right)
  \exp\left(-\frac{5 
      e^{-4}(x-K/2)^4}{12 T} \right) dx \left(1 + \mathcal{O}(T) \right) \\
  &= \frac{3^{1/4} \Gamma(1/4)}{5^{1/4} 2 \sqrt{2} \pi e^2}
  \exp\left(-\frac{1}{T}\right) \frac{1}{T^{3/4}} \left(1 + \mathcal{O}(T)
  \right),
\end{align*}
where we used
\begin{equation*}
  \int_{-\infty}^\infty \exp(-\alpha x^4) dx = \frac{\Gamma(1/4)}{2
    \alpha^{1/4}}, \quad \alpha > 0.
\end{equation*}
Thus, we arrive at~\eqref{eq:density-expansion-b}.

\section{Large Deviations approach}
\label{sec:large-devi-appr}

Our main tool here are novel marginal density expansions in small-noise regime
\cite{DFJVpartI}. This was used in order to compute the large-strike behavior
of implied volatility in the correlated Stein--Stein model;
\cite{Stein91,GuSt}.\footnote{Similar investigations have recently been
  conducted in the Heston model; \cite{He, Gu12} and the references therein.}

In fact, the technical assumptions of \cite{DFJVpartI} were satisfied in the
analysis of the Stein--Stein model whereas in the (seemingly) trivial case of
two IID Black-Scholes assets, the technical assumptions of \cite{DFJVpartI}
are indeed violated for a critical strike $K=K^{\ast }$. The necessity of this
condition is then highlighted by the fact, as was seen in the previous
section,
\begin{equation*}
\left\{ \frac{\partial ^{2}}{\partial K^{2}}C_{B}\left( K,T\right) \right\}
_{K=K^{\ast }} \nsim \text{\textrm{(const)}} \exp \left( -\frac{\Lambda \left(
K^{\ast }\right) }{T}\right) \frac{1}{T^{1/2}}.
\end{equation*}
The computation of $K^{\ast }$ can be achieved either via a geometric
construction borrowed from Riemannian geometry, which relies on the Weingarten
map, or by some (fairly) elementary analysis of a system of Hamiltonian
ODEs. In fact, the Hamiltonian point of view extends naturally when one
introduces correlation, local and even stochastic volatility. Explicit answers
then depend on the analytical tractability of these (boundary value) ODE
problems. (Of course, the numerical solution of such problems is well-known.)

In the following, we review \cite{DFJVpartI}. Consider a $d$-dimensional
diffusion $\left( \mathrm{X}_{t}^{\varepsilon }\right) _{t\geq 0}$ given by
the stochastic differential equation
\begin{equation}
d\mathrm{X}_{t}^{\varepsilon }=b\left( \varepsilon ,\mathrm{X}%
_{t}^{\varepsilon }\right) dt+\varepsilon \sigma \left( \mathrm{X}%
_{t}^{\varepsilon }\right) dW_{t},\quad \text{with }\mathrm{X}%
_{0}^{\varepsilon }=\mathrm{x}_{0}^{\varepsilon }\in \mathbb{R}^{d},
\label{SDEXeps}
\end{equation}%
and where $W=(W^{1},\dots ,W^{m}$) is an $m$-dimensional Brownian motion.
Unless otherwise stated, we assume $b:[0,1)\times \mathbb{R}^{d}\rightarrow
\mathbb{R}^{d},$ $\sigma =\left( \sigma _{1},\dots ,\sigma _{m}\right) :%
\mathbb{R}^{d}\rightarrow \mathrm{Lin}\left( \mathbb{R}^{m},\mathbb{R}%
^{d}\right) $ and $\mathrm{x}_{0}^{\cdot }:[0,1)\rightarrow \mathbb{R}^{d}$
to be smooth, bounded with bounded derivatives of all orders. Set $\sigma
_{0}=b\left( 0,\cdot \right) $ and assume that, for every multi-index $\alpha
$, the drift vector fields $b\left( \varepsilon ,\cdot \right) $ converges
to $\sigma _{0}$ in the sense\footnote{%
If (\ref{SDEXeps}) is understood in Stratonovich sense, so that $dW$ is
replaced by $\circ dW$, the drift vector field $b\left( \varepsilon ,\cdot
\right) $ is changed to $\tilde{b}\left( \varepsilon ,\cdot \right) =b\left(
\varepsilon ,\cdot \right) -\left( \varepsilon ^{2}/2\right)
\sum_{i=1}^{m}\sigma _{i}\cdot \partial \sigma _{i}$. In particular, $\sigma
_{0}$ is also the limit of $\tilde{b}\left( \varepsilon ,\cdot \right) $ in
the sense of (\ref{bepsTob}) .}%
\begin{equation}
\partial _{x}^{\alpha }b\left( \varepsilon ,\cdot \right) \rightarrow
\partial _{x}^{\alpha }b\left( 0,\cdot \right) =\partial _{x}^{\alpha
}\sigma _{0}\left( \cdot \right) \text{ uniformly on compacts as }%
\varepsilon \downarrow 0\text{.}  \label{bepsTob}
\end{equation}%
We shall also assume that%
\begin{equation}
\partial _{\varepsilon }b\left( \varepsilon ,\cdot \right) \rightarrow
\partial _{\varepsilon }b\left( 0,\cdot \right) \text{ \ uniformly on
compacts as }\varepsilon \downarrow 0  \label{bepsC1}
\end{equation}%
and
\begin{equation}
\mathrm{x}_{0}^{\varepsilon }=\mathrm{x}_{0}+\varepsilon \mathrm{\hat{x}}%
_{0}+o\left( \varepsilon \right) \text{ as }\varepsilon \downarrow 0\text{. }
\label{x0eps_ass}
\end{equation}

\begin{theorem}
\label{thm:MainThm} \textbf{(Small noise)} Let $\left( \mathrm{X}%
^{\varepsilon }\right) $ be the solution process to%
\begin{equation*}
d\mathrm{X}_{t}^{\varepsilon }=b\left( \varepsilon ,\mathrm{X}%
_{t}^{\varepsilon }\right) dt+\varepsilon \sigma \left( \mathrm{X}%
_{t}^{\varepsilon }\right) dW_{t},\quad \text{with }\mathrm{X}%
_{0}^{\varepsilon }=\mathrm{x}_{0}^{\varepsilon }\in \mathbb{R}^{d}.
\end{equation*}%
Assume $b\left( \varepsilon ,\cdot \right) \rightarrow \sigma _{0}\left( \cdot
\right) $ in the sense of (\ref{bepsTob}), (\ref{bepsC1}), and $%
\mathrm{X}_{0}^{\varepsilon }\equiv \mathrm{x}_{0}^{\varepsilon }\rightarrow
\mathrm{x}_{0}$ as $\varepsilon \rightarrow 0$ in the sense of (\ref%
{x0eps_ass}). Assume non-degeneracy of $\sigma$ in the sense that
$\sigma. \sigma^T$ is strictly positive definite everywhere in space.
\footnote{This may be relaxed to a weak Hoermander condition with an explicit
  controllability condition.}  Fix $\mathrm{y}\in
\mathbb{R}^{l},\,N_{\mathrm{y}}:=\left( \mathrm{y},\cdot \right) $ and let
$\mathcal{K}_{\mathrm{y}}$ be the the space of all $%
\mathrm{h}\in H$, the Cameron-Martin space of absolutely continuous paths with
derivatives in $L^2\left([0,T], \mathbb{R}^m \right)$, s.t. the solution to%
\begin{equation*}
d\phi _{t}^{\mathrm{h}}=\sigma _{0}\left( \phi _{t}^{\mathrm{h}}\right)
dt+\sum_{i=1}^{m}\sigma _{i}\left( \phi _{t}^{\mathrm{h}}\right) d\mathrm{h}%
_{t}^{i},\,\,\phi _{0}^{\mathrm{h}}=\mathrm{x}_{0}\in \mathbb{R}^{d}
\end{equation*}%
satisfies $\phi _{T}^{\mathrm{h}}\in N_{\mathrm{y}}$. In a neighborhood of $%
\mathrm{y}$, assume smoothness of\footnote{%
If $\#\mathcal{K}_{\mathrm{y}}^{\min }=1$ smoothness of the energy can be
shown and need not be assumed; \cite{DFJVpartI}. Note also that in our
application to tail asymptotics, with $\theta $-scaling, $\theta \in \left\{
1,2\right\} $, the energy must be linear resp. quadratic (by scaling) and
hence smooth.}
\begin{equation*}
\Lambda \left( \mathrm{y}\right) =\inf \left\{ \frac{1}{2}\Vert \mathrm{h}%
\Vert _{H}^{2}:\mathrm{h}\in \mathcal{K}_{\mathrm{y}}\right\} .
\end{equation*}%
Assume also (i) there are only finitely many minimizers, i.e. $\mathcal{K}_{%
\mathrm{y}}^{\min }<\infty $ where
\begin{equation*}
\mathcal{K}_{\mathrm{y}}^{\min }:=\left\{ \mathrm{h}_{0}\in \mathcal{K}_{%
\mathrm{y}}:\frac{1}{2}\Vert \mathrm{h}_{0}\Vert _{H}^{2}=\Lambda \left(
\mathrm{y}\right) \right\} ;
\end{equation*}%
(ii) $\mathrm{x}_{0}$ is non-focal for $N_{\mathrm{y}}$ in the sense of
\cite{DFJVpartI}. (We shall review below how to check this.)  Then
there exists $c_{0}=c_{0}\left(
  \mathrm{x}_{0},\mathrm{y},T\right) >0$ such that
\begin{equation*}
\mathrm{Y}_{T}^{\varepsilon }=\Pi _{l}\mathrm{X}_{T}^{\varepsilon }=\left(
X_{T}^{\varepsilon ,1},\dots ,X_{T}^{\varepsilon ,l}\right) ,\,\,\,\,1\leq
l\leq d,
\end{equation*}%
admits a density with expansion%
\begin{equation*}
f_{\varepsilon }\left( \mathrm{y},T\right) =e^{-\frac{\Lambda \left( \mathrm{%
y}\right) }{\varepsilon ^{2}}}e^{\,\frac{\max \left\{ \Lambda ^{\prime
}\left( \mathrm{y}\right) \cdot \,\mathrm{\hat{Y}}_{T}\left( \mathrm{h}%
_{0}\right) :\mathrm{h}_{0}\in \mathcal{K}_{\mathrm{y}}^{\min }\right\} }{%
\varepsilon }}\varepsilon ^{-l}\left( c_{0}+O\left( \varepsilon \right)
\right) \text{ as }\varepsilon \downarrow 0,
\end{equation*}%
where $\Lambda'$ denotes the gradient of $\Lambda$.

Here $\mathrm{\hat{Y}=\hat{Y}}\left( \mathrm{h}_{0}\right) =\left( \hat{Y}%
^{1},\dots ,\hat{Y}^{l}\right) $ is the projection, $\mathrm{\hat{Y}=}\Pi
_{l}\mathrm{\hat{X}}$, of the solution to the following (ordinary)
differential equation
\begin{eqnarray}
d\mathrm{\hat{X}}_{t} &=&\Big(\partial _{x}b\left( 0,\phi _{t}^{\mathrm{h}%
_{0}}\left( \mathrm{x}_{0}\right) \right) +\partial _{x}\sigma (\phi _{t}^{%
\mathrm{h}_{0}}\left( \mathrm{x}_{0}\right) )\mathrm{\dot{h}}_{0}\left(
t\right) \Big)\mathrm{\hat{X}}_{t}dt+\partial _{\varepsilon }b\left( 0,\phi
_{t}^{\mathrm{h}_{0}}\left( \mathrm{x}_{0}\right) \right) dt,
\label{eq:SDEYHat} \\
\,\,\,\,\,\mathrm{\hat{X}}_{0} &=&\mathrm{\hat{x}}_{0}.  \notag
\end{eqnarray}
\end{theorem}

\begin{remark}[Localization]
The assumptions on the coefficients $b,\sigma $ in theorem \ref{thm:MainThm}
(smooth, bounded with bounded derivatives of all orders) are typical in this
context (cf. Ben\ Arous \cite{Benarous, Benarous2} for instance) but rarely
met in practical examples from finance. This difficulty can be resolved by a
suitable localization. For instance, as detailed in \cite{DFJVpartI}, an
estimate of the form
\begin{equation}
\lim_{R\rightarrow \infty }\lim \sup_{\varepsilon \rightarrow 0}\varepsilon
^{2}\log \mathbb{P}\left[ \tau _{R}\leq T\right] =-\infty .  \label{Assloc}
\end{equation}%
with $\tau _{R}:=\inf \left\{ t\in \left[ 0,T\right] :\sup_{s\in \left[ 0,t%
\right] }\left\vert \mathrm{X}_{s}^{\varepsilon }\right\vert \geq R\right\} $
will allow to bypass the boundedness assumptions.
\end{remark}

\subsection{Short time asymptotics}

The reduction of \textit{short time expansions} to small noise expansions by
Brownian scaling is classical. In the present context, we have the following
statement, taken from \cite[Sec. 2.1]{DFJVpartI}.

\begin{corollary}
  \textbf{(Short time)} \label{CorShortTime}Consider $d\mathrm{X}_{t}=b\left(
    \mathrm{X}_{t}\right) dt+\sigma \left( \mathrm{X}_{t}\right) dW$, started
  at $\mathrm{X}_{0}=\mathrm{x}_{0}\in \mathbb{R}^{d}$, with $C^{\infty }$%
  -bounded vector fields which are non-degenerate in the sense that
  $\sigma. \sigma^T$ is strictly positive definite everywhere in space. Fix
  $\mathrm{y}\in \mathbb{R}^{l},\,N_{\mathrm{y}}:=\left( \mathrm{y},\cdot
  \right) $ and assume (i),(ii) as in theorem \ref{thm:MainThm}. Let $%
  f\left( t,\cdot \right) =f\left( t,\mathrm{y}\right) $ be the density of $%
  \mathrm{Y}_{t}=\left( \mathrm{X}_{t}^{1},\dots ,\mathrm{X}_{t}^{l}\right) $.
  Then
\begin{equation*}
f\left( t,\mathrm{y}\right) \sim \text{(const)}\frac{1}{t^{l/2}}\exp \left( -%
\frac{d^{2}\left( \mathrm{x}_{0},\mathrm{y}\right) }{2t}\right) \text{ as }%
t\downarrow 0
\end{equation*}%
where $d\left( \mathrm{x}_{0},\mathrm{y}\right) $ is the sub-Riemannian
distance, based on $\left( \sigma _{1},\dots ,\sigma _{m}\right) $, from the
point $\mathrm{x}_{0}$ to the affine subspace $N_{\mathrm{y}}$.
\end{corollary}

\subsection{Computational aspects}
\label{sec:comp-aspects}

We present here the \textit{mechanics} of the actual computations, in the
spirit of the Pontryagin maximum principle (e.g. \cite{SS}).
For details we refer to \cite{DFJVpartI}.

\begin{itemize}
\item \textbf{The Hamiltonian. }Based on the SDE (\ref{SDEXeps}), with
diffusion vector fields $\sigma _{1},\dots ,\sigma _{m}$ and drift vector
field $\sigma _{0}$ (in the $\varepsilon \rightarrow 0$ limit) we define the
\textit{Hamiltonian}%
\begin{align*}
\mathcal{H}\left( \mathrm{x},\mathrm{p}\right) &\coloneqq \left\langle
  \mathrm{p},\sigma _{0}\left( \mathrm{x}\right) \right\rangle +\frac{1}{2}
\sum_{i=1}^{m}\left\langle \mathrm{p},\sigma _{i}\left( \mathrm{x}\right)
\right\rangle ^{2} \\ 
&=\left\langle \mathrm{p},\sigma _{0}\left( \mathrm{x}\right) \right\rangle
+\frac{1}{2}\left\langle \mathrm{p},\left( \sigma \sigma ^{T}\right) \left(
\mathrm{x}\right) \mathrm{p}\right\rangle .
\end{align*}
Remark the driving Brownian motions $W^{1},\dots ,W^{m}$ were assumed to be
independent. Many stochastic models, notably in finance, are written in terms
of correlated Brownian motions, i.e. with a non-trivial correlation matrix
$\Omega =\left( \omega ^{i,j}:1\leq i,j\leq m\right) $, where $d\,\left\langle
  W^{i},W^{j}\right\rangle _{t}=\omega ^{i,j}dt$. The Hamiltonian then
becomes%
\begin{equation}
\mathcal{H}\left( \mathrm{x},\mathrm{p}\right) =\left\langle \mathrm{p}%
,\sigma _{0}\left( \mathrm{x}\right) \right\rangle +\frac{1}{2}\left\langle
\mathrm{p},\left( \sigma \Omega \sigma ^{T}\right) \left( \mathrm{x}\right)
\mathrm{p}\right\rangle .  \label{HamiltonianWithCorrel}
\end{equation}

\item \textbf{The Hamiltonian ODEs.} The following system of ordinary
differential equations,
\begin{equation}
\left(
\begin{array}{c}
\mathrm{\dot{x}(t)} \\
\mathrm{\dot{p}(t)}%
\end{array}%
\right) =\left(
\begin{array}{c}
\partial _{\mathrm{p}}\mathcal{H}\left( \mathrm{x}\left( t\right) ,\mathrm{p}%
\left( t\right) \right) \\
-\partial _{\mathrm{x}}\mathcal{H}\left( \mathrm{x}\left( t\right) ,\mathrm{p%
}\left( t\right) \right)%
\end{array}%
\right) ,  \label{HamiltonianODEs}
\end{equation}%
gives rise to a solution flow, denoted by $\mathrm{H}_{t\leftarrow 0}$, so
that%
\begin{equation*}
\mathrm{H}_{t\leftarrow 0}\left( \mathrm{x}_{0},\mathrm{p}_{0}\right)
\end{equation*}%
is the unique solution to the above ODE with initial data $\left( \mathrm{x}
  _{0},\mathrm{p}_{0}\right) $. Our standing (regularity) assumption are more
than enough to guarantee uniqueness and local ODE\ existence. As in \cite[
p.37]{Bismut}, the vector\ field $\left( \partial _{\mathrm{p}}\mathcal{H}
  ,-\partial _{\mathrm{x}}\mathcal{H}\right) $ is complete, i.e., one has global
existence. It can be useful to start the flow backwards with time-$T$ terminal
data, say $\left( \mathrm{x}_{T},\mathrm{p}_{T}\right) $; we then write
\begin{equation*}
\mathrm{H}_{t\leftarrow T}\left( \mathrm{x}_{T},\mathrm{p}_{T}\right)
\end{equation*}%
for the unique solution to (\ref{HamiltonianODEs}) with given time-$T$
terminal data. Of course,%
\begin{equation*}
\mathrm{H}_{t\leftarrow T}\left( \mathrm{H}_{T\leftarrow 0}\left( \mathrm{x}%
_{0},\mathrm{p}_{0}\right) \right) =\mathrm{H}_{t\leftarrow 0}\left( \mathrm{%
x}_{0},\mathrm{p}_{0}\right) .
\end{equation*}

\item \textbf{\ Solving the Hamiltonian ODEs as boundary value problem. }%
Given the target manifold $N_{\mathrm{a}}=\left( \mathrm{a},\cdot \right) $,
the analysis in \cite{DFJVpartI} requires solving the Hamiltonian ODEs (\ref%
{HamiltonianODEs}) with mixed initial -, terminal - and transversality
conditions,%
\begin{align}
\mathrm{x}\left( 0\right) &=\mathrm{x}_{0}\in \mathbb{R}^{d},  \notag \\
\mathrm{x}\left( T\right) &=\left( \mathrm{y},\cdot \right) \in \mathbb{R}%
^{l}\mathbb{\oplus R}^{d-l}, \label{ITTcond} \\
\mathrm{p}\left( T\right) &=\left( \cdot ,0\right) \in \mathbb{R}^{l}%
\mathbb{\oplus R}^{d-l}. \notag 
\end{align}%
Note that this is a $2d$-dimensional system of ordinary differential
equations, subject to $d+l+\left( d-l\right) =2d$ conditions. In general,
boundary problems for such ODEs may have more than one, exactly one or no
solution. In the present setting, there will always be one or more than one
solution. After all, we know by \cite{DFJVpartI} that there exists at least
one minimizing control $\mathrm{h}_{0}$ and that can be reconstructed via the
solution of the Hamiltonian ODEs, as explained in the following step.

\item \textbf{Finding the minimizing controls.} The Hamiltonian ODEs, as
boundary value problem, are effectively first order conditions (for
minimality) and thus yield \textit{candidates} for the minimizing control $%
\mathrm{h}_{0}=\mathrm{h}_{0}\left( \cdot \right) $, given by%
\begin{equation}
\mathrm{\dot{h}}_{0}=\left(
\begin{array}{c}
\,\left\langle \sigma _{1}\left( \mathrm{x}\left( \cdot \right) \right) ,%
\mathrm{p}\left( \cdot \right) \right\rangle \\
\dots \\
\,\,\left\langle \sigma _{m}\left( \mathrm{x}\left( \cdot \right) \right) ,%
\mathrm{p}\left( \cdot \right) \right\rangle%
\end{array}
\right) .  \label{Formula_h0}
\end{equation}
Each such candidate is indeed admissible in the sense $\mathrm{h}_{0}\in
\mathcal{K}_{\mathrm{a}}$ but may fail to be a minimizer. We thus compute the
energy $\left\Vert \mathrm{h}_{0}\right\Vert _{H}^{2} = \mathcal{H}(x_0, p_0)$
for each candidate and identify those (``$\mathrm{h}_{0}\in
\mathcal{K}_{\mathrm{a}}^{\min }$'') with minimal energy. The procedure via
Hamiltonian flows also yields a unique $\mathrm{p}_{0}=\mathrm{p}_{0}\left(
  \mathrm{h}_{0}\right) $. If $\sigma_0 = 0$ -- as in our case -- the energy
is equal to $\mathcal{H}(x_0, p_0)$, otherwise the formula is slightly more
complicated.

\item \textbf{Checking non-focality. }By definition \cite{DFJVpartI}, $%
\mathrm{x}_{0}$ is \textbf{non-focal} for $N=\left( \mathrm{y},\cdot \right)
$ along $\mathrm{h}_{0}\in \mathcal{K}_{\mathrm{a}}^{\min }$ in the sense
that, with $\left( \mathrm{x}_{T}\mathrm{,p}_{T}\right) :=\mathrm{H}%
_{T\leftarrow 0}\left( \mathrm{x}_{0}\mathrm{,p}_{0}\left( \mathrm{h}%
_{0}\right) \right) \in \mathcal{T}^{\ast }\mathbb{R}^{d}$,
\begin{equation*}
\partial _{\left( \mathfrak{z},\mathfrak{q}\right) }|_{\left( \mathfrak{z},%
\mathfrak{q}\right) \mathfrak{=}\left( 0,0\right) }\pi \mathrm{H}%
_{0\leftarrow T}\left( \mathrm{x}_{T}+\left(
\begin{array}{c}
0 \\
\mathfrak{z}%
\end{array}%
\right) ,\mathrm{p}_{T}+\left( \mathfrak{q},0\right) \right)
\end{equation*}%
is non-degenerate (as $d\times d$ matrix; here we think of $\left( \mathfrak{%
    z},\mathfrak{q}\right) \in \mathbb{R}^{d-l}\times \mathbb{R}^{l}\cong
\mathbb{R}^{d}$ and recall that $\pi $ denotes the projection from $\mathcal{%
  T}^{\ast }\mathbb{R}^{d}$ onto $\mathbb{R}^{d}$; in coordinates $\pi \left(
  \mathrm{x},\mathrm{p}\right) =\mathrm{x}$). Note that in the point-point
setting, $\mathrm{x}_{T}=\mathrm{y}$ is fixed and only perturbations of the
arrival "velocity" $\mathrm{p}_{T}$ - without restrictions, i.e. without
transversality condition - are considered. Non-degeneracy of the resulting map
should then be called \textbf{non-conjugacy} (between two points; here: $%
\mathrm{x}_{T}$ and $\mathrm{x}_{0}$). In the absence of the drift vector
field $\sigma_0$, this is consistent with the usual meaning of non-conjugacy;
after identifying tangent- and cotangent-space $\partial
_{\mathfrak{q}}|_{\mathfrak{q=}0}\pi \mathrm{H}_{0\leftarrow T}$ is precisely
the differential of the exponential map.

\item \textbf{The explicit marginal density expansion. }We then have
\begin{equation*}
\text{ }f^{\varepsilon }\left( \mathrm{y},T\right) =e^{-c_{1}/\varepsilon
^{2}}e^{c_{2}/\varepsilon }\varepsilon ^{-l}\left( c_{0}+O\left( \varepsilon
\right) \right) \text{ as }\varepsilon \downarrow 0.
\end{equation*}%
with $c_{1}=\Lambda \left( \mathrm{y}\right) $. The second-order exponential
constant $c_{2}$ then requires the solution of a finitely many ( $\#\mathcal{%
K}_{\mathrm{a}}^{\min }<\infty $) auxiliary ODEs, cf. theorem \ref%
{thm:MainThm}.
\end{itemize}

\section{Analysis of the Black--Scholes basket}
\label{sec:analys-black-schol}

For a general multi-dimensional Black-Scholes model, we have a Hamiltonian
$$\mathcal{H}(x,p) =  \frac{1}{2} \left\langle p, (\sigma(x) \Omega
  \sigma(x)^T) p \right\rangle, $$ with $\sigma(x) = (\sigma^1 x^1, \ldots,
\sigma^m x^m)$. While the corresponding Hamiltonian ODEs can be solved in
closed form, the boundary conditions lead to systems of non-linear equations,
which we cannot solve explicitly any more. While numerical solutions are, of
course, possible, we restrict ourselves to the extremely simple setting of
Section~\ref{sec:comp-based-saddle}, in order to keep maximal tractability.

Consequently, we have the Hamiltonian $\mathcal{H}(x,p) = \frac{1}{2}
\left( (\sigma x^1 p^1)^2 + (\sigma x^2 p^2)^2\right)$. The solutions of the
Hamiltonian ODEs started at $(x_0,p_0)$ satisfy
\begin{equation}
  \label{eq:hamiltonian-ode-solution}
  H_{t \leftarrow 0}(x_0,p_0) =
  \begin{pmatrix}
    x^1_0 e^{\sigma^2 x^1_0 p^1_0 t}\\
    x^2_0 e^{\sigma^2 x^2_0 p^2_0 t}\\
    p^1_0 e^{-\sigma^2 x^1_0 p^1_0 t}\\
    p^2_0 e^{-\sigma^2 x^2_0 p^2_0 t}
  \end{pmatrix},
\end{equation}
which can be easily seen from the observation that $\mathcal{H}$ is constant
along solutions of the Hamiltonian ODEs together with symmetry between
$(x^1,p^1)$ and $(x^2,p^2)$. This immediately implies that the inverse flow is
given by
\begin{equation}
  \label{eq:hamiltonian-inverse-flow}
  H_{0 \leftarrow t}(x_t,p_t) =
  \begin{pmatrix}
    x^1_t e^{-\sigma^2 x^1_t p^1_t t}\\
    x^2_t e^{-\sigma^2 x^2_t p^2_t t}\\
    p^1_t e^{\sigma^2 x^1_t p^1_t t}\\
    p^2_t e^{\sigma^2 x^2_t p^2_t t}
  \end{pmatrix}.
\end{equation}
Now we introduce the boundary conditions. Note that, contrary to
Theorem~\ref{thm:MainThm}, we now project to the linear subspace
$\{x:x^1+x^2=K\}$. Thus, the terminal condition on $x$ translates into $x^1_T +
x^2_T = K$ -- we need to end at the target manifold --, whereas the
transversality condition translates to $p_T$ being orthogonal to the target
manifold. Evaluating these conditions at $T = 1$, we get
\begin{gather*}
  x^1_0 = S^1_0 = 1,\\
  x^2_0 = S^2_0 = 1,\\
  x^1_1 + x^2_1 = K,\\
  p^1_1 - p^2_1 = 0.
\end{gather*}
It is a pleasant exercise to check that solving for $x^1_1=:x$ and $x^2_x = K-x$
then leads exactly to the first order condition \eqref{eq:LaplaceFirstOrderCond}
encountered in Section \ref{sec:comp-based-saddle}. With identical arguments,
 assuming $K \le 2e$ from here on (and disregarding the case $K>2e$ where closed form computations are not available),
we find that the optimal configuration must satisfy $x_1^\ast
= (K/2, K/2)$. Inserting this value into the first two components
of~\eqref{eq:hamiltonian-ode-solution}, we obtain the equation
\begin{equation*}
  \frac{K}{2} = e^{\sigma^2 p^i_0} \iff p^i_0 = \log\left( \frac{K}{2}
  \right)/\sigma^2, \quad i=1,2.
\end{equation*}
This implies that $p^\ast_1 = \left( \frac{2}{\sigma^2K} \log(K/2),
  \frac{2}{\sigma^2K} \log(K/2) \right)$.
Moreover, we see that the minimizing control satisfies
\begin{equation}\label{eq:h_0}
  \dot{h}_0(t) =
  \begin{pmatrix}
    \sigma x^1(t) p^1(t)\\
    \sigma x^2(t) p^2(t)
  \end{pmatrix}
  =
  \begin{pmatrix}
    \sigma p^1_0\\
    \sigma p^2_0
  \end{pmatrix}
  =
  \begin{pmatrix}
    \f{\log(K/2)}{\sigma}\\
    \f{\log(K/2)}{\sigma}
  \end{pmatrix},
\end{equation}
see~\eqref{Formula_h0}, implying that the minimal energy is given by
\begin{equation}
  \label{eq:minimal-energy-example}
  \Lambda(K) = \half \norm{h_0}_H^2 = \f{\log(K/2)^2}{\sigma^2} =
  \mathcal{H}(x_0,p_0).
\end{equation}

Regarding focality, we have to check that the matrix:
\begin{equation}
  \label{eq:matrix-derivatives}
  M(x_1,p_1) :=
  \begin{pmatrix}
    \left.\frac{\partial}{\partial \epsilon} \right|_{\epsilon = 0} H^1_{0
      \leftarrow 1}(x_1 + \epsilon (1,-1), p_1) & \left.\frac{\partial}{\partial
        \eta} \right|_{\eta = 0} H^1_{0 \leftarrow 1}(x_1, p_1 + \eta (1,1)) \\
    \left.\frac{\partial}{\partial \epsilon} \right|_{\epsilon = 0} H^2_{0
      \leftarrow 1}(x_1 + \epsilon (1,-1), p_1) &
    \left.\frac{\partial}{\partial
        \eta} \right|_{\eta = 0} H^2_{0 \leftarrow 1}(x_1, p_1 + \eta (1,1))
  \end{pmatrix}
\end{equation}
is non-degenerate when evaluated at the optimal configuration $(x_1^\ast,
p_1^\ast)$. A simple calculation shows that
\begin{equation*}
  M(x_1, p_1) =
  \begin{pmatrix}
    e^{-\sigma^2 x_1^1 p_1^1} - x^1_1 p_1^1 \sigma^2 e^{-\sigma^2 x_1^1 p_1^1}
    & -\sigma^2 (x_1^1)^2 e^{-\sigma^2 x_1^1p_1^1} \\
    -e^{-\sigma^2 x_1^2 p_1^2} + x^2_1 p_1^2 \sigma^2 e^{-\sigma^2 x_1^2 p_1^2}
    & -\sigma^2 (x_1^2)^2 e^{-\sigma^2 x_1^2 p_1^2}
  \end{pmatrix},
\end{equation*}
implying that
\begin{equation*}
  M(x_1^\ast, p_1^\ast) =
  \begin{pmatrix}
    \frac{2}{K}(1-\log(K/2)) & -\frac{\sigma^2K}{2}\\
    \frac{2}{K}(-1+\log(K/2)) & -\frac{\sigma^2K}{2}
  \end{pmatrix},
\end{equation*}
and we can conclude that
\begin{equation*}
  \det M(x_1^\ast, p_1^\ast) = 2 \sigma^2 \left(\log(K/2) - 1 \right),
\end{equation*}
which is zero if and only if $K = 2e$. We summarize the results of this
calculation as follows:
\begin{itemize}
\item In the generic case $K \neq 2e$, the non-focality condition of
  Theorem~\ref{thm:MainThm} holds true, and we obtain (from Corollary~\ref{CorShortTime}) 
  the following (short time) density expansion of $B_T = \exp(\sigma W^1_T) + \exp (\sigma W^2_T)$, 
  expansion
  $$
    K \mapsto \exp \left( -\frac{\Lambda \left( K\right) }{T}\right)
      \frac{1}{\sqrt{T}} \left( c_0+O\left( T\right) \right)
  $$
  When specialized to unit volatility, we recover precisely ~\eqref{eq:density-expansion-a}.
\item For $K = 2e$, the initial stock price is focal for the minimizing
  configuration, so the non-focality condition of Theorem~\ref{thm:MainThm}
  fails. And indeed, we {\it want} it to fail for the actual expansion in this case, 
  namely ~\eqref{eq:density-expansion-b}, is not at all of the generic form predicted by
  our theorem.
\end{itemize}

\begin{remark}
  It is immediate to use this analysis to deal also with the case of non-unit
  (but identical) spots $S^1_0 = S^2_0$ by scaling the Black-Scholes dynamics
  accordingly, i.e., by replacing $K$ with $K/S_0^1$. Hence, in this case
  focality happens when $\log\left(\f{K}{2S_0^1}\right) = 1$, i.e., when $K =
  2 S_0^1 e$.
\end{remark}

\begin{remark} The question arises if the critical (``focal") case $K = 2e$, with atypical algebraic factor $T^{-3/4}$
cf. ~\eqref{eq:density-expansion-b}, can also be recovered by a general theorem. Related results in 
\cite{Mo75} and also \cite{TW} suggests that this may indeed be the case but would require substantial additional
work.
\end{remark}

\section{Extensions: correlation, local and stochastic vol}
\label{sec:extens-corr-local}

\subsection{Analysis of the Black--Scholes basket, small noise}
In section \ref{sec:analys-black-schol} we analyzed the density of a simple
Black--Scholes basket with dynamics 
$$
dB_t = S^1_t \sigma dW^1_t + S^2_t \sigma dW^2_t.
$$ 
As explained in Section~\ref{sec:large-devi-appr} the analysis is really based
on a small noise (small vol) expansion 
of
$$
dB^\epsilon_t = S^{1,\epsilon}_t \sigma \epsilon dW^1_t + S^{2,\epsilon}_t
\sigma \epsilon dW^2_t, 
$$ 
run til time $T=1$. Consider now a situation with small rates, also of order
$\epsilon$. In other words, 
$$
dS^{i,\epsilon}_t = r S^{i,\epsilon}_t \epsilon dt + S^{i,\epsilon}_t \sigma
\epsilon dW^i, 
$$  
and then $B^\epsilon_t = S^{1,\epsilon}_t + S^{2,\epsilon}_t$ as before. We
still assume $S_0^i = 1$. A look at Theorem~\ref{thm:MainThm} (now we cannot
use Corollary~\ref{CorShortTime}) reveals that the entire leading order
computation remains unchanged (at least at unit time and with trivial changes
otherwise). The resulting (now: small noise) density expansion of
$B^\epsilon_T|_{T=1}$ is more involved and takes the form
\begin{equation}
  \label{eq:small-noise-example}
  K \mapsto \exp \left( -\frac{\Lambda \left( K\right) }{\epsilon^2}\right) 
  \exp \left( \f{2r\log(K/2)}{\sigma^2 \log(2) \epsilon}\right)
  \frac{1}{\epsilon} \left( c_0 +O(\epsilon) \right). 
\end{equation}
Here $\Lambda \left( K\right)$ is given in closed form,
cf.~\eqref{eq:minimal-energy-example}, so that $\Lambda'\left( K\right) =
\frac{2 \log(K/2)}{\sigma^2 K}$ is also explicitly known. Furthermore, under
similar restrictions on $K$ as before, $h_0$
is (still) given by~\eqref{eq:h_0}, so that
\begin{equation*}
  \phi^{h_0}_t =
  \begin{pmatrix}
    (K/2)^t\\
    (K/2)^t
  \end{pmatrix}.
\end{equation*}
Thus, the ODE for $\hat{X}$ (see Theorem~\ref{thm:MainThm}) is given by
\begin{equation*}
 \f{d\hat{X}_t}{dt} = \log(K/2) \hat{X}_t + r
 \begin{pmatrix}
   (K/2)^t\\
   (K/2)^t
 \end{pmatrix}, \quad \hat{X}_0 = \hat{x}_0 = 0,
\end{equation*}
which has the solution
\begin{equation*}
  \hat{X}_t^i = r\left(1- \left(\half \right)^t \right) \f{K^t}{\log 2},
\end{equation*}
implying that $\hat{Y}_1 = \hat{X}_1^1 + \hat{X}_1^2 = rK/\log(2)$. Thus, the
second exponential term has the form given above.

\subsection{Basket analysis under local, stochastic vol etc.}

One can immediately write down the Hamiltonian associated to, say two, or $d>2$ assets, each of which is governed by 
local vol dynamics or stochastic vol, based on additional factors. In general, however, one will be stuck with the 
analysis of the resulting boundary value problem for the Hamiltonian ODEs; numerical (e.g. shooting) methods will 
have to be used. In  some models, including the Stein--Stein model, we believe (due to the analysis carried out 
in \cite{{DFJVpartII}}) that, in special cases, closed form answers are possible but we will not pursue this here. 
Instead, we continue with a few more computation in the Black--Scholes case for $d$ assets.

\subsection{Multi-variate Black-Scholes models}
\label{sec:multi-variate-black}

In the multi-variate case $d > 2$ of a general, $d$-dimensional Black Scholes
model with correlation matrix $(\rho_{ij})$, the Hamiltonian has the form
\begin{equation*}
  \mathcal{H}(x,p) = \half \sum_{i,j=1}^d \rho_{ij} \sigma^i p^i x^i \sigma^j
  x^j p^j.
\end{equation*}
Thus, the Hamiltonian ODEs have the form
\begin{align*}
  \dot{x}^l &= \sigma^l x^l \sum_{i=1}^d \rho_{li} \sigma^i p^i x^i, \quad
  i=1, \ldots, d\\
  \dot{p}^l &= -\sigma^l p^l \sum_{i=1}^d \rho_{li} \sigma^i p^i x^i, \quad
  i=1, \ldots, d.
\end{align*}
Consequently, it is again easy to see that $\f{\pa}{\pa t} x^l(t) p^l(t) = 0$,
implying that $x^l(t) p^l(t) = x^l_0 p^l_0$. The Hamiltonian flow has the form
\begin{equation}
  \label{eq:hamiltonian-flow-d-cor}
  H_{t \leftarrow 0}(x_0,p_0) =
  \begin{pmatrix}
    \left( x^l_0 \exp\left[ \sigma^l \left( \sum_{i=1}^d \rho_{li} \sigma^i
          p^i_0 x^i_0 \right) t \right] \right)_{l=1}^d\\
    \left( p^l_0 \exp\left[ -\sigma^l \left( \sum_{i=1}^d \rho_{li} \sigma^i
          p^i_0 x^i_0 \right) t \right] \right)_{l=1}^d
  \end{pmatrix}.
\end{equation}
Using again that $p^l(t)x^l(t) = p^l(0) x^l(0)$ for any $l$, we obtain the
inverse Hamiltonian flow
\begin{equation}
  \label{eq:hamiltonian-invflow-d-cor}
  H_{0 \leftarrow t}(x_t,p_t) =
  \begin{pmatrix}
    \left( x^l_t \exp\left[ -\sigma^l \left( \sum_{i=1}^d \rho_{li} \sigma^i
          p^i_t x^i_t \right) t \right] \right)_{l=1}^d\\
    \left( p^l_t \exp\left[ \sigma^l \left( \sum_{i=1}^d \rho_{li} \sigma^i
          p^i_t x^i_t \right) t \right] \right)_{l=1}^d
  \end{pmatrix}.
\end{equation}
The boundary conditions -- at $T=1$ -- are now given by
\begin{subequations}
  \label{eq:boundary-d-dim}
  \begin{gather}
    x_0 = S_0\label{eq:boundary-d-dim-a}\\
    \sum_{l=1}^d x^l(1) = K\label{eq:boundary-d-dim-b}\\
    p^1(1) = p^2(1) = \cdots = p^d(1).\label{eq:boundary-d-dim-c}
  \end{gather}
\end{subequations}
Indeed, the transversality condition~\eqref{eq:boundary-d-dim-c} says that the
final momentum $p(1)$ is orthogonal to the surface $\Set{\sum_{l=1}^d y^l =
  K}$, whose tangent space is spanned by the collection of vectors
$\mathbf{e}_1 - \mathbf{e}_l$, $l=2, \ldots, d$, with $\mathbf{e}_1, \ldots,
\mathbf{e}_d$ the standard basis of $\R^d$. The
equations~\eqref{eq:boundary-d-dim} are certainly not difficult to solve
numerically, but an explicit solution is not available, neither in the general
case nor in the case of $d$ uncorrelated assets.
\begin{remark}
  The main point of this calculation is that while explicit solutions are no
  longer possible in a general Black-Scholes model, the
  phenomenon~\eqref{eq:density-expansion} potentially appears in all
  Black-Scholes models. Moreover, we stress that the non-focality conditions
  are easily checked numerically.
\end{remark}
\begin{remark}
  \label{rem:asian-discrete}
  Note that the discretely monitored Asian option can be considered as a
  special case of a basket option on correlated assets. Indeed, let us
  consider an option on
  \begin{equation*}
    \f{1}{N} \sum_{i=1}^N S_{t_i}, \text{ with (for simplicity) } t_i = i
    \dt, \quad i=1, \ldots, N.
  \end{equation*}
  For each individual $i \in \set{1, \ldots, N}$ we have, for fixed $\dt>0$,
  the equality in law 
  \begin{equation*}
    S_{t_i} = S_0 e^{\sigma B_{i \dt} - \half \sigma^2 i \dt} = S_0
    e^{\sigma^i W^i_{\dt} - \half (\sigma^i)^2 \dt}
  \end{equation*}
  for $\sigma^i \coloneqq \sqrt{i} \sigma$ and $W^i_{\dt} \coloneqq
  B_{i\dt}/\sqrt{i}$. In law, the vector $\left(W^1_{\dt}, \ldots, W^N_{\dt}
  \right)$ corresponds to the marginal distribution of an $N$-dimensional
  Brownian motion at time $\dt$ with correlation $\rho_{ij} =
  \f{\min(i,j)}{\sqrt{ij}}$, $1 \le i,\, j\le N$. Thus, the Asian option
  corresponds to an option on the basket with $S_0^i \equiv S_0$, $\sigma^i$
  as above and a correlation matrix $\rho_{ij}$ with maturity $\dt$. Moreover,
  the asymptotic expansion of the price of the Asian option as $\dt \to 0$
  corresponds to the short-time asymptotics of the basket.
\end{remark}
\begin{remark}
  \label{rem:asian-continuous}
  A small-noise asymptotic expansion of the continuous Asian option on
  $\int_0^T S_t dt$ is also possible by the techniques of
  Section~\ref{sec:large-devi-appr} (with ellipticity conditions replaced by
  weak H\"{o}rmander conditions). Essentially, this is equivalent to letting
  $N \to \infty$ in Remark~\ref{rem:asian-discrete} -- but more direct.
\end{remark}

As in the two-dimensional case, the boundary conditions can be solved
explicitly in the fully symmetric case, when $\sigma^l \equiv \sigma$ and,
say, $S^l_0 \equiv 1$. For suitable $K$ the optimal configuration is
\begin{gather*}
  x_0^\ast = (1, \ldots, 1)^T, \quad x_1^\ast = (K/d, \ldots, K/d)^T\\
  p_0^\ast = \left( \f{\log(K/d)}{\sigma^2}, \ldots, \f{\log(K/d)}{\sigma^2}
  \right)^T, \quad p_1^\ast = \left( \f{d}{\sigma^2 K} \log(K/d), \ldots,
    \f{d}{\sigma^2 K} \log(K/d) \right)^T.
\end{gather*}
Introducing 
\begin{equation*}
  \mathfrak{q} = \epsilon_1
  \begin{pmatrix}
    1\\
    1\\
    \vdots\\
    1
  \end{pmatrix},\quad
  \mathfrak{z} =
  \begin{pmatrix}
    \epsilon_2 + \cdots + \epsilon_d\\
    -\epsilon_2\\
    \vdots\\
    -\epsilon_d
  \end{pmatrix},
\end{equation*}
we obtain (for the case of $d$ uncorrelated assets)
\begin{align*}
  M(x_1, p_1) &\coloneqq \left.\pa_{(\mathfrak{z},
      \mathfrak{q})}\right|_{(\mathfrak{z}, \mathfrak{q}) = 0} \pi H_{0
    \leftarrow 1}(x_1 + \mathfrak{z}, p_1 + \mathfrak{q}) \\
  &=
  \begin{pmatrix}
    a_1 & \mathbf{b} \\
    \mathbf{a} & G
  \end{pmatrix},
\end{align*}
where $\mathbf{a} = (a_2, \ldots, a_d)^T \in \R^{(d-1)\times 1}$, $\mathbf{b}
= b (1, \ldots, 1) \in \R^{1\times (d-1)}$, $G = \operatorname{diag}(g_2,
\ldots, g_d) \in \R^{(d-1)\times(d-1)}$ with 
\begin{gather*}
  a_l = -(\sigma^l)^2 (x_1^l)^2 e^{-(\sigma^l)^2 p_1^l x_1^l}, \quad l=1,
  \ldots d,\\
  b = \left[ 1 - (\sigma^1)^2 x_1^1 p_1^1 \right] e^{-(\sigma^1)^2 p_1^1
    x_1^1}, \\
  g_l = - \left[ 1 - (\sigma^l)^2 x_1^l p_1^l \right] e^{-(\sigma^l)^2 p_1^l
      x_1^l}, \quad l=2,\ldots d.
\end{gather*}
In the symmetric case, we can evaluate $M$ at the optimal configuration and
obtain
\begin{equation*}
  M(x_1^\ast, p_1^\ast) =
  \begin{pmatrix}
    -\sigma^2 \f{K}{d} & \left[1 - \log(K/d) \right] \f{d}{K} & \cdots &
    \left[1 - \log(K/d) \right] \f{d}{K} \\
    -\sigma^2 \f{K}{d} & -\left[1 - \log(K/d) \right] \f{d}{K} & \cdots & 0 \\ 
    \vdots & \vdots & \ddots & \vdots\\
    -\sigma^2 \f{K}{d} & 0 & \cdots & -\left[1 - \log(K/d) \right] \f{d}{K}
  \end{pmatrix},
\end{equation*}
whose determinant can be seen to be
\begin{equation*}
  \det M(x_1^\ast, p_1^\ast) = (-1)^{d} \sigma^2 K \left[ \left( 1 - \log(K/d)
    \right) \f{d}{K} \right]^{d-1}. 
\end{equation*}
Thus, the non-focality condition fails if and only if $K = d e$. Moreover, we
obtain the energy
\begin{equation*}
  \Lambda(K) = \mathcal{H}(x^\ast_0, p^\ast_0) = \f{d}{2}
  \f{\log(K/d)^2}{\sigma^2}.
\end{equation*}

\section{A geometric approach to focality}
\label{sec:geom-appr-focal}

In this final section we take a more geometrical look at the non-focality
condition appearing in Section~\ref{sec:comp-aspects}. Consider the Black
Scholes model
\begin{equation*}
  dS^i_t = \sigma^iS^i_t dW^i_t,\quad \left\langle dW^i, dW^j \right\rangle_t
  = \rho_{i,j} dt. 
\end{equation*}
We change parameters $\mathbf{S} \to \mathbf{y} \to \mathbf{x}$,
by
\begin{equation*}
  y^i := \frac{\log\left( \frac{S^i}{S^i_0} \right)}{\sigma^i}, \quad x^i =
  L_{ip} y^p,\quad i=1, \ldots, d,
\end{equation*}
where $\boldsymbol{\rho}$ denotes the correlation matrix of $\mathbf{W}$ and
$\boldsymbol{\rho} = L L^T$ its Cholesky factorization. Obviously, $S^i =
S_0^i e^{\sigma_i y^i}$. In terms of the $\mathbf{x}$-coordinates we have
\begin{gather*}
  x^i = x^i(\mathbf{F}) = L_{ip} \log\left(S^p/S_0^p\right)/\sigma^p,\\
  S^i = S^i(\mathbf{x}) = S_0^i e^{\sigma_i L^{ip} x^p}.
\end{gather*}
The advantage of using the chart $\mathbf{x}$ is that the corresponding
Riemannian metric tensor is the usual Euclidean metric tensor. Thus, we
simply have
\begin{equation*}
  d(\mathbf{S}_0, \mathbf{S}) = |\mathbf{x}_0 - \mathbf{x}|
\end{equation*}
and the geodesics are straight lines as seen from the $\mathbf{x}$-chart. Note
furthermore that $\mathbf{S} = \mathbf{S}_0$ is transformed to $\mathbf{x} =
\mathbf{0}$.


The payoff function of the option is given by $\left(\sum w_i S^i_T - K
\right)^+$. We normalize $w_i \equiv 1$ and $T \equiv 1$. The strike surface
$F = \left\{\mathbf{S}\in \mathbb{R}^d_+ \left| \sum_{i=1}^d S^i =
    K\right. \right\}$, which is (a sub-set of) a hyperplane in $\mathbf{S}$
coordinates is, however, transformed to a much more complicated submanifold in
$\mathbf{x}$ coordinates. Re-phrasing the equation $\sum_i S^i = K$ in
$\mathbf{y}$-coordinates and solving for $y^d$ gives
\begin{equation*}
  y^d = \log\left[ \left(K - \sum_{i=1}^{d-1} S^i_0 e^{\sigma^i\sum_{p=1}^d L^{ip} x^p}
    \right) / S_0^d \right] / \sigma^d,
\end{equation*}
with $(L^{ij}) = (L_{ij})^{-1}$,
which implies -- using that $L$ and $L^{-1}$ are lower-triangular matrices --
\begin{equation*}
  L^{dd} x^d = \log\left[ \left(K - \sum_{i=1}^{d-1} S_0^i e^{\sigma^i
        \sum_{p=1}^i L^{ip} x^p} \right) / S_0^d \right] / \sigma^d -
  \sum_{k=1}^{d-1} L^{dk} x^k.
\end{equation*}
For sake of clarity, let us introduce the notation $\mathbf{q} = (q^1, \ldots,
q^{d-1}) := (x^1, \ldots, x^{d-1})$. A parametrization of the strike
surface $F$ is then given by the map $\varphi: U \subset \mathbb{R}^{d-1} \to
\mathbb{R}^d$ with 
\begin{equation*}
  U := \left\{\mathbf{q} \in \mathbb{R}^{d-1} \left| \sum_{i=1}^{d-1} S_0^i
    e^{\sigma^i \sum_{p=1}^i L^{ip} q^p} < K \right. \right\},
\end{equation*}
and
\begin{equation*}
  \varphi(\mathbf{q}) := \left(\mathbf{q},\, \frac{1}{L^{dd}} \left\{
      \log\left[ \left(K - \sum_{i=1}^{d-1} S_0^i e^{\sigma^i 
        \sum_{p=1}^i L^{ip} q^p} \right) / S_0^d \right] / \sigma^d -
  \sum_{k=1}^{d-1} L^{dk} q^k \right\} \right).
\end{equation*}
Note that by the change of coordinates, we are implicitly assuming that $S^i >
0$ for all $i$. Moreover, the standard basis $\mathbf{e}_1(\mathbf{p}),
\ldots, \mathbf{e}_{d-1}(\mathbf{p})$ of the tangent space $T_{\mathbf{p}}F$
to $F$ at $\mathbf{p} = \varphi(\mathbf{q})$ is given by the columns of the
Jacobi matrix of $\varphi$ evaluated at $\mathbf{q}$, more precisely we have
\begin{equation*}
  \mathbf{e}_i(\mathbf{p}) = \left( (\delta_i^j)_{j=1}^{d-1},\, -\frac{1}{L^{dd}}
  \left[ \frac{1}{\sigma^d} \frac{\sum_{j=i}^{d-1} \sigma^j L^{ji} S_0^j
      e^{\sigma^j \sum_{r=1}^j L^{jr} q^r}}{K - \sum_{j=1}^{d-1} S_0^j
      e^{\sigma^j \sum_{r=1}^j L^{jr} q^r}} + L^{di}\right]\right)
\end{equation*}
for $i=1, \ldots, d-1$ and $\mathbf{p} = \varphi(\mathbf{q})$. Consequently,
the normal vector field $N$ to $S$ at $\mathbf{p} = \varphi(\mathbf{q})$ is
given by
\begin{equation*}
  N(\mathbf{p}) = \alpha(\mathbf{p}) \left( \left( \frac{1}{L^{dd}}
      \left[ \frac{1}{\sigma^d} \frac{\sum_{j=i}^{d-1} \sigma^j L^{ji} S_0^j
          e^{\sigma^j \sum_{r=1}^j L^{jr} q^r}}{K - \sum_{j=1}^{d-1} S_0^j
          e^{\sigma^j \sum_{r=1}^j L^{jr} q^r}} + L^{di}\right]
    \right)_{i=1}^{d-1},\, 1 \right) = N\circ \varphi(\mathbf{q}),
\end{equation*}
where $\alpha$ is a normalization factor guaranteeing that
$|N(\mathbf{p})| = 1$, i.e.,
\begin{equation*}
  \alpha(\mathbf{p}) = \left(1 + \sum_{i=1}^{d-1} \frac{1}{(L^{dd})^2}
      \left[ \frac{1}{\sigma^d} \frac{\sum_{j=i}^{d-1} \sigma^j L^{ji} S_0^j
          e^{\sigma^j \sum_{r=1}^j L^{jr} q^r}}{K - \sum_{j=1}^{d-1} S_0^j
          e^{\sigma^j \sum_{r=1}^j L^{jr} q^r}} + L^{di}\right]^2 \right)^{-1/2}.
\end{equation*}

The \emph{Weingarten map} or \emph{shape operator} $L_{\mathbf{p}}:
T_{\mathbf{p}}F \to T_{\mathbf{p}}F$ is defined by
\begin{equation*}
  L_{\mathbf{p}}\left( d\varphi_{\varphi^{-1}(\mathbf{p})}( \mathbf{v} )
  \right) = - d (N \circ \varphi) (\varphi^{-1}(\mathbf{p})) \cdot \mathbf{v},
\end{equation*}
$\mathbf{v} \in \mathbb{R}^{d-1} = T_{\varphi^{-1}(\mathbf{p})} U$,
see~\cite{doCarmo}. In other words, for $\varphi(\mathbf{q}) = \mathbf{p}$, we
interpret $N$ as a map in $\mathbf{q}$ and $-L_{\mathbf{p}}$ is the
directional derivative of that map. We study the Weingarten map since it gives
us the curvature of the surface $F$. Indeed, the eigenvalues $k_1(\mathbf{p}),
\ldots, k_{d-1}(\mathbf{p})$ of the linear map $L_{\mathbf{p}}:
T_{\mathbf{p}}F \to T_{\mathbf{p}}F$ are called \emph{principal curvatures} of
$F$. Then the \emph{focal points} of $F$ at $\mathbf{p}$ are given by
\begin{equation*}
  \{\mathbf{p} + \frac{1}{k_i(\mathbf{p})} N(\mathbf{p}) | 1 \le i \le d-1
    \text{ such that } k_i(\mathbf{p}) \neq 0\}.
\end{equation*}

In order to compute the eigenvalues of the shape operator, we need to compute
the representation of $L_{\mathbf{p}}$ in the standard basis
$\left(\mathbf{e}_1(\mathbf{p}),\ldots, \mathbf{e}_{d-1}(\mathbf{p})
\right)$. Let us denote this matrix by $\overline{L}(\mathbf{p})$, then we
obviously have
\begin{equation*}
  \overline{L}(\varphi(\mathbf{q}))_{ij} = -\langle \frac{\partial}{\partial
    q^j} (N \circ \varphi)(\mathbf{q}), \mathbf{e}_i(\varphi(\mathbf{q}))
  \rangle,\quad i,j=1, \ldots, d-1. 
\end{equation*}
The principal curvatures $k_1(\mathbf{p}), \ldots, k_{d-1}(\mathbf{p})$ are,
thus, the eigenvalues of the $(d-1)$-dimensional matrix
$\overline{L}(\mathbf{p})$. 

Since the calculations become too complicated in the general case, we now
again concentrate on the case of \emph{two uncorrelated} assets, i.e., $d=2$
and $\boldsymbol{\rho} = L = I_2$. In this case, we have
\begin{gather*}
  \mathbf{e}_1(\mathbf{p}) = \left(1,\, -\frac{\sigma^1}{\sigma^2}
    \frac{S_0^1 e^{\sigma^1 q^1}}{K-S_0^1 e^{\sigma^1 q^1}} \right), \\
  N(\varphi(\mathbf{q})) = \frac{1}{\sqrt{(\sigma^1)^2 (S_0^1)^2 e^{2 \sigma^1
        q^1} 
    + (\sigma^2)^2 \left(K - S_0^1 e^{\sigma^1 q^1} \right)^2}} \left(
  \sigma^1 S_0^1 e^{\sigma^1 q^1},\, \sigma^2 \left(K - S_0^1 e^{\sigma^1
      q^1} \right) \right).
\end{gather*}
Thus, the Weingarten map is given by
\begin{equation*}
  L_{\mathbf{p}}(v \mathbf{e}_1(\mathbf{p})) = v \kappa(\mathbf{p})
  \mathbf{e}_1(\mathbf{p}), 
\end{equation*}
where for $\mathbf{q} = (q^1) \in \mathbb{R}$
\begin{equation*}
  \kappa(\varphi(\mathbf{q})) = k_1(\varphi(\mathbf{q})) = \frac{K (\sigma^1)^2
    (\sigma^2)^2 S_0^1 e^{\sigma^1
    q^1} \left( S_0^1 e^{\sigma^1 q^1} - K \right)}{\left[(\sigma^1)^2
    (S_0^1)^2 e^{2 \sigma^1 q^1} + (\sigma^2)^2 \left(S_0^1 e^{\sigma^1 q^1} -
      K\right)^2 \right]^{3/2}}
\end{equation*}
is the \emph{curvature} of the curve $F$ in $\mathbb{R}^2$. We see that
$\kappa = 0$ if and only if $K = S_0^1 e^{\sigma^1 q^1}$, i.e., at the
boundary of the surface $F$. Otherwise, $\kappa$ is negative.

Here, both components of $N(\mathbf{p})$ are positive on
$F$. Consequently, for any $\mathbf{p} = \varphi(\mathbf{q}) \in S$ there is
precisely one focal point $\mathfrak{f} = \mathfrak{f}(\mathbf{p}) \in
\mathbb{R}^2$, which is given by
\begin{align*}
  \mathfrak{f}^1 &= q^1 + \frac{S_0^1 e^{\sigma^1 q^1} \left[2 (\sigma^2)^2 K -
      ((\sigma^1)^2 + (\sigma^2)^2) S_0^1 e^{\sigma^1 q^1} \right] - (\sigma^2)^2
    K^2}{\sigma^1 (\sigma^2)^2 K\left(K - S_0^1 e^{\sigma^1 q^1}\right)},\\
  \mathfrak{f}^2 &= \frac{1}{\sigma^2} \log\left( \frac{K- S_0^1 e^{\sigma^1
        q^1}}{S_0^2} \right) + 2 \frac{\sigma^2}{(\sigma^1)^2} -
  \frac{\sigma^2 K
    e^{-\sigma^1 q^1}}{(\sigma^1)^2S_0^1} - \frac{((\sigma^1)^2+(\sigma^2)^2) S_0^1
  e^{\sigma^1 q^1}}{(\sigma^1)^2 \sigma^2 K}.
\end{align*}
Denoting $\mathbf{p} = (x^1, x^2)$ and re-introducing the short-cut notation
$S^i = S_0^i e^{\sigma^i x^i}$,
$i=1,2$, (noting that $S^1+S^2=K$) we can express $\mathfrak{f}$ as 
\begin{align*}
  \mathfrak{f}^1 &= x^1 + \frac{S^1 \left[ 2 (\sigma^2)^2 K -
      ((\sigma^1)^2+(\sigma^2)^2) S^1 \right] - (\sigma^2)^2 K^2}{\sigma^1
    (\sigma^2)^2 K S^2},\\
  \mathfrak{f}^2 &= x^2 + \frac{S^1 \left[ 2 (\sigma^2)^2 K -
      ((\sigma^1)^2+(\sigma^2)^2) S^1 \right] - (\sigma^2)^2 K^2}{(\sigma^1)^2
    \sigma^2 K S^1}.
\end{align*}

In the current setting, let $\mathbf{q}^\ast$ be the optimal configuration in
$\mathbf{q}$-coordinates, i.e., the point on $F$ with smallest Euclidean
norm. Then the non-focality condition of Theorem~\ref{thm:MainThm} is
satisfied, if $0$ is not a focal point to
$\varphi(\mathbf{q}^\ast)$, see the discussion in the proof of
\cite[Prop.~6]{DFJVpartI}.
\begin{remark}
  As both components of the normal vector $N$ are non-negative on $F$ and the
  curvature $\kappa$ is negative, $\mathbf{0}$ can only be a focal point if
  $F$ has a non-empty intersection with the positive quadrant. Inserting into
  the parametrization of $F$, we see that this can only be the case if $K >
  S_0^1 + S_0^2$. In other words: if the option is in the money, then the
  non-focality condition is always satisfied (in the two-dimensional,
  uncorrelated case).
\end{remark}
Let us again use the parameters of
Section~\ref{sec:comp-based-saddle}, i.e., $S_0^1=S_0^2=1$, $\sigma^1=\sigma^2
= \sigma$. Then we consider $\mathbf{S}^\ast = (K/2, K/2)$, which translates
into $\mathbf{x}^\ast = \left(\frac{\log(K/2)}{\sigma},
  \frac{\log(K/2)}{\sigma} \right)$. Inserting into the formulas for the focal
points, we obtain
\begin{equation*}
  \mathfrak{f}^1(\mathbf{x}^\ast) =  \mathfrak{f}^2(\mathbf{x}^\ast) =
  \frac{\log\left(\frac{K}{2}\right) - 1}{\sigma}.
\end{equation*}
So, $0$ is focal to the optimal configuration, if and only if
\begin{equation*}
  K = 2e,
\end{equation*}
and we recover, once more, the results of Section~\ref{sec:comp-based-saddle}
and Section~\ref{sec:analys-black-schol} -- recall that $\mathbf{S}_0$
corresponds to $0$ in $\mathbf{x}$-coordinates.

\begin{figure}[!htb]
  \centering
  \subfloat[Optimal configuration (out of money regime)]{
    \includegraphics[width=0.45\textwidth]{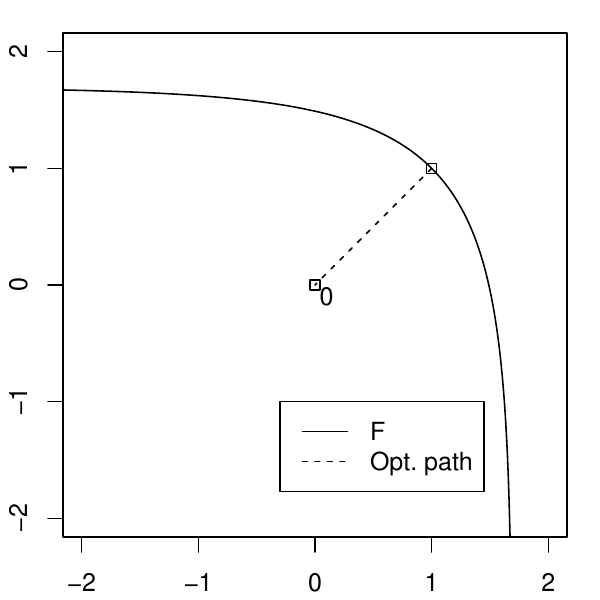}
  }
  \qquad
  \subfloat[Focal points]{
    \includegraphics[width=0.45\textwidth]{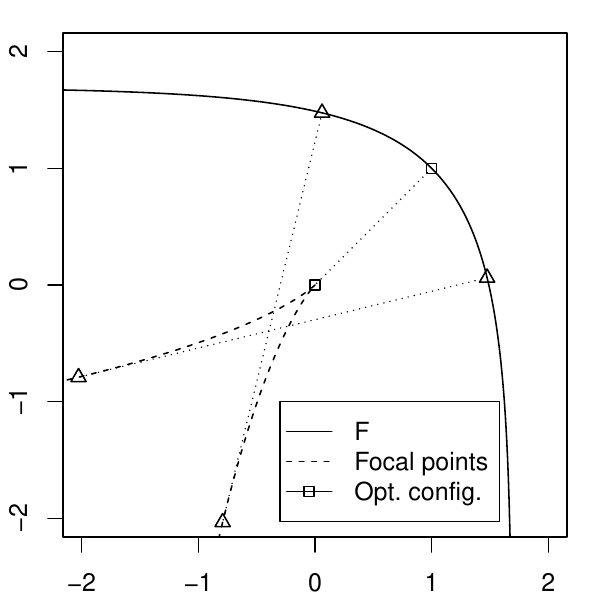}
  }
  \caption{Optimal configuration and focal points for two independent assets
    with $\sigma^1 = \sigma^2 = 1$, $\mathbf{S}_0 = (1,1)$, $K = 2e$.\\
    (A) The dashed line depicts the optimal path between the spot price
    $\mathbf{S}_0$ ($0$ in the $\mathbf{q}$-chart) and the optimal
    configuration.\\ 
    (B) Dotted lines connect some selected points on the manifold $F$ with the
    corresponding focal points. Points marked with a triangle visualize the
    construction of the focal points. We see that $0$ is, indeed, focal to the
  optimal configuration.}
  \label{fig:two-uncorrelated-1}
\end{figure}

\begin{figure}[!htb]
  \centering
  \subfloat[Optimal configuration (in the money regime)]{
    \includegraphics[width=0.45\textwidth]{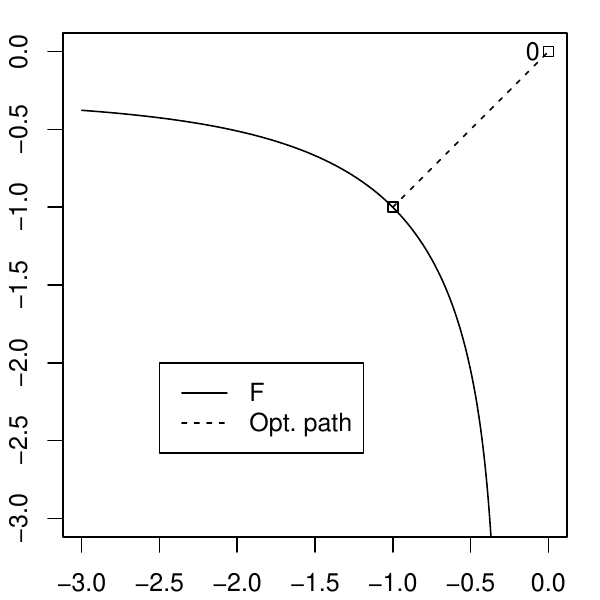}
  }
  \qquad
  \subfloat[Focal points]{
    \includegraphics[width=0.45\textwidth]{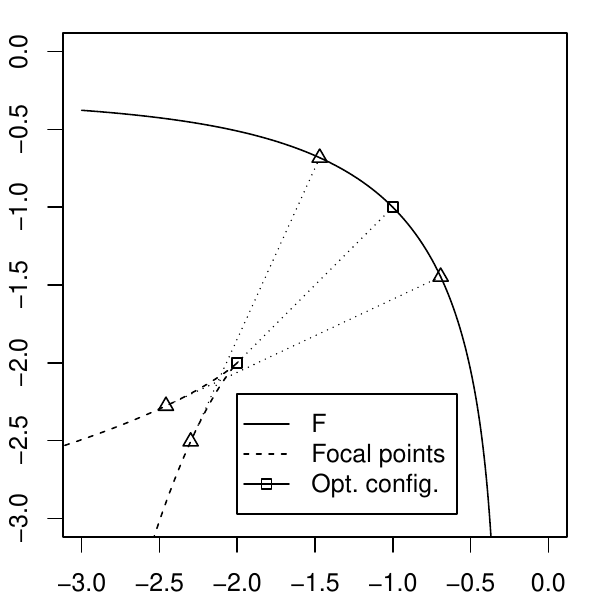}
  }
  \caption{Optimal configuration and focal points for two independent assets
    with $\sigma^1 = \sigma^2 = 1$, $\mathbf{S}_0 = (1,1)$, $K = 2/e$.\\
    (A) The dashed line depicts the optimal path between the spot price
    $\mathbf{S}_0$ ($0$ in the $\mathbf{q}$-chart) and the optimal
    configuration.\\ 
    (B) Dotted lines connect some selected points on the manifold $F$ with the
    corresponding focal points. Points marked with a triangle visualize the
    construction of the focal points. This example illustrates the fact that
    the non-focality condition always holds when the basket option is in the
    money.}
  \label{fig:two-uncorrelated-2}
\end{figure}

In Figure~\ref{fig:two-uncorrelated-1} and~\ref{fig:two-uncorrelated-2} the
focal points are visualized for two different configurations of two
uncorrelated baskets. We plot the surface $F$ as a submanifold of
$\mathbb{R}^2$. We have seen above that for any $\mathbf{p} \in F$ there is
precisely one focal point $\mathfrak{f}(\mathbf{p})$. Hence, we additionally
plot the surface $\{\mathfrak{f}(\mathbf{p}) | \mathbf{p} \in F\}$ -- more
precisely, part of this surface. In Figure~\ref{fig:two-uncorrelated-1} we
show the case constructed above where the non-focality condition is
violated. In Figure~\ref{fig:two-uncorrelated-2} the option is ITM. As
explained above, in the ITM case the manifold $F$ does not intersect the
positive quadrant, implying that the non-focality condition is satisfied.

\end{document}